\documentclass[12pt]{amsart}

\author{Stefano Decio}
\address{Department of Mathematical Sciences, Norwegian University of Science and Technology, 7491 Trondheim, Norway}
\email{stefano.decio@ntnu.no}

\usepackage{amsmath}
\usepackage{amssymb}
\usepackage{amsfonts}
\usepackage{enumerate}
\usepackage{mathrsfs}

\newcommand{\R}{{\mathbf R}}

\newcommand{\ld}{{\lambda}}

\newtheorem{theorem}{Theorem}
\newtheorem{lemma}{Lemma}
\newtheorem{proposition}{Proposition}
\theoremstyle{remark}
\newtheorem*{remark}{Remark}
\theoremstyle{definition}

\DeclareMathOperator{\diver}{div}
\begin{document}

\begin{abstract}
   We consider linear combinations of eigenfunctions of the Laplace-Beltrami operator on a compact Riemannian manifold $(M,g)$ and investigate a density property of their zero sets. More precisely, let $f=\sum_{k=1}^m a_k \phi_{\ld_{j_k}}$, where $-\Delta_g\phi_{\ld}=\ld\phi_{\ld}$. Denoting by $Z_f$ the zero-set of $f$, we show that for any $x\in M$, $dist(x,Z_f)\leq C(m)\ld_{j_1}^{-1/2}$. The proof is based on a new integral Harnack-type estimate for positive solutions of higher order elliptic PDEs.
   
\smallskip
\noindent \textbf{Keywords.} Laplace eigenfunctions; nodal sets; higher order PDEs.

\smallskip
\noindent \textbf{Mathematics Subject Classification (2010).} Primary 35J30; Secondary 35R01, 58J05, 58J50.
\end{abstract}

\title{Density of zero sets for sums of eigenfunctions}
\maketitle

\section{Introduction}

Let $(M,g)$ be a $\mathcal{C}^{\infty}$-smooth compact Riemannian manifold of dimension $d$. Our object of study are the eigenfunctions of (minus) the Laplace-Beltrami operator $\Delta_g$ on $M$, that is solutions of the equation
\begin{align}
    \Delta_g\phi_{\ld}+\ld\phi_{\ld}=0.
\end{align}
It is well known that the eigenvalues are positive, form a discrete set and can be arranged in a sequence $\{\ld_n\}$ that tends to infinity; moreover, the normalized eigenfunctions form an orthonormal basis for $L^2(M)$.\\

An interesting fact about eigenfunctions is that their zero sets $Z_{\phi_\ld}=\{x\in M: \phi_{\ld}(x)=0\}$ are $\frac{Constant}{\sqrt{\ld}}$ dense in $M$, meaning that for any $y\in M$ there holds:
$$dist(y,Z_{\phi_\ld})\leq \frac{C(M,g)}{\sqrt{\ld}}$$
and examples show that this is sharp; note that this implies an upper bound of $C/\sqrt{\ld}$ for the radius of the largest ball inscribed in a nodal domain of $\phi_{\ld}$. There are several different proofs of this fact: a simple argument based on domain monotonicity for eigenvalues is contained in \cite{B}, but, as pointed out by the referee, it was probably already known to Polya and Szeg{\"o}, or maybe even before by Weyl. We will highlight a different argument in the next section. \\
The inverse problem is in this case more interesting and still the subject of research: can one inscribe a ball  of radius $\sim\ld^{-1/2}$ in a nodal domain? If one interprets the question as being about a fixed domain $\Omega$, and asking whether one can inscribe a ball of radius $\sim \ld_{1}(\Omega)^{-1/2}$ (where $\ld_{1}$ is here the first Dirichlet eigenvalue) in $\Omega$, then the answer is in general false in dimension greater than $3$, due to the possibility of having thin spikes in the domain that do not have much impact on $\ld_1$ but affect the inner radius considerably. However, in the celebrated paper \cite{Lb} Lieb showed that one can inscribe 'most' of such a ball inside $\Omega$; improvements are contained in \cite{GM}, \cite{LS}. In two dimensions the result is true for a simply connected domain, see \cite{Ma} and \cite{Ha}. The situation is different if $\Omega$ is not just any domain, but a nodal domain of a global eigenfunction on $M$; in this case, essentially due to controlled local asymmetry of signs of the eigenfunction, there exist some lower bounds on the inner radius, see \cite{M}, \cite{GM}. Except in dimension $2$, the best available bound is not the conjectured optimal one of $C/\sqrt{\ld}$; see also \cite{G} for an improvement in case the metric is real analytic.\\

Now, one could ask if a similar density property holds for the zero set of a linear combination of eigenfunctions corresponding to different eigenvalues, where the lowest eigenvalue is large enough (such a sum has to vanish somewhere, by orthogonality to constants). After all, if we think of eigenfunctions as vast generalizations of homogeneous trigonometric polynomials, it is only natural to consider a sum of them. 
Some thought suggests that perhaps this is false if too many eigenfunctions are considered; the problem is already nontrivial on the torus.
\subsection*{Trigonometric Polynomials}
In \cite{KO} it is proved that a real multivariate trigonometric polynomial 
$$f(x)=\sum_{\nu\in S}c(\nu)\exp(2\pi i \langle x,\nu\rangle)$$
where $0\notin S\subset \mathbf{Z}^d$, $S=-S$, has a zero in every ball of radius $$R(S)=\sum_{\nu\in S}\frac{1}{8|\nu|}.$$
Note that we can view $f$ as a linear combination of eigenfunctions on the flat torus $\mathbf{T}^d$, with eigenvalue $\ld=|\nu|^2$; then the result is in general not optimal due to the possible multiplicity of $\ld$: $R(S)$ could be large if $|\nu|^2=\ld$ for every $\nu\in S$, while in this case $f$ is a single eigenfunction for which we know the right density scale should be $\sim \ld^{-1/2}$.\\
On the other hand, the theorem is sharp in the following sense: the authors construct in one dimension a trigonometric polynomial with spectrum equal to, say, $([-2N,N]\cup[N,2N])\cap \mathbf{Z}$ which is strictly positive on some fixed sub-interval of the unit circle with size independent of $N$.\\
We also note a related result of Eremenko-Novikov in \cite{EN}, answering a conjecture of Benjamin Logan: the result essentially says that if a function on the real line has a large spectral gap (in the sense of Fourier transform support) at the origin, there have to be many sign changes of the function. We thank the referee for pointing out this circle of ideas and the reference \cite{EN}.\\

The proof in \cite{KO} is very specific to the special case of the torus; we are interested in whether a similar result holds true for a general compact smooth manifold. Our main result offers a partial answer. \\
\textbf{Notational warning:} To avoid the proliferation of indices, we use a somewhat non-standard notation: in this paper, $\lambda_1$ is \textit{not} the first eigenvalue of the Laplace operator, but merely the lowest eigenvalue in the linear combination we consider. \\
The dependence of the constants on the dimension $d$ is pervasive throughout the paper and will not be explicitly indicated. 

\begin{theorem}
\label{main}
Let $(M,g)$ be a $\mathcal{C}^{\infty}$-smooth compact Riemannian manifoldas, and let $f=\sum_{j=1}^{m}a_j\phi_{{\ld}_j}$ where the $a_j$'s are real numbers. Suppose $\ld_1\leq \ld_2\leq...\leq \ld_m $. Then, if $\ld_1>C_1(M,g,m)$, $f$ has a zero in every geodesic ball of radius $$\frac{C(M,g,m)}{\sqrt{\ld_1}},$$ where
\begin{align*}
    C(M,g,m)= \begin{cases}
    C(M,g)m^{\frac{d+2}{2}}\quad \text{if} \quad d\geq 3\\
    C(M,g)m\log^2m\quad \text{if}\quad d=1,2.      
    \end{cases}
\end{align*}
\end{theorem}

\begin{remark}
 The example in \cite{KO} shows that one cannot hope for better than a linear dependence of the density radius on $m$. The super-linear dependence on $m$ of $C(M,g,m)$ when $d\geq 3$ is likely to be an artifact of the proof. 
 
 For a general smooth metric we do not know the growth (in $m$) of $C_1(M,g,m)$, since in our method it depends on $2m-1$ derivatives of the metric $g$. It is conceivable that for a real analytic metric one could show that $C_1(M,g,m)$ is linear in $m$, which would be optimal; however, the combinatorics of keeping track of the derivatives does not seem to be as straightforward as one might hope. In the real analytic case we have some estimates for $C_1(M,g,m)$ but we decided to not include them as the they are probably far from being optimal. For $m$ bounded independently of $\ld_1$ Theorem \ref{main} is sharp. 
\end{remark}

\subsection*{Related work} Results of Fang-Hua Lin in \cite{Li} (see also Jerison-Lebeau \cite{JL} and Donnelly \cite{D}) show that a linear combination of eigenfunctions cannot 'vanish much more' than the highest frequency eigenfunction in the sum, in various precise senses. In the opposite direction, Stefan Steinerberger (see \cite{S}, \cite{S2}) proved some bounds from below for the Hausdorff measure of the zero set of functions orthogonal to the first few eigenfunctions. He sets his result in the context of a tentative generalization of Sturm-Hurwitz theory to higher dimensions; the present work could be seen in the same light. For a modern account of the 1-dimensional theory see \cite{BH}. Note however that density and lower measure bounds can be quite different characteristics of nodal sets, as the next example shows.

\subsection*{Isolated zeros}
Consider the sphere $\mathbf{S}^d$ with $d\geq 2$. Take two zonal spherical harmonics, rotationally invariant with respect to $x_0\in \mathbf{S}^d$, $Z_{k}$ and $Z_{n}$, corresponding to eigenvalues $k(k+d-1)$ and $n(n+d-1)$. It is known that they both have a peak at $x_0$ with $d_k:=Z_k(x_0)\sim k^{\frac{d-1}{2}}$, $d_n:=Z_n(x_0)\sim n^{\frac{d-1}{2}}$. Suppose that $k<n$. Let $$f=Z_k-\frac{d_k}{d_n}Z_n;$$
Theorem \ref{main} shows that $f$ has zeros everywhere at density scale $\sim k^{-1}$; on the other hand, $x_0$ is an isolated zero for $f$. By changing the coefficient in front of $Z_n$ one can get a component of the zero set around $x_0$ with arbitrarily small diameter. This does not happen with a single eigenfunction and appears to be an essential difficulty in obtaining lower measure bounds for the zero set, even in two dimensions. Theorem \ref{main} gives an upper bound for the inner radius of a nodal domain of a linear combination of eigenfunctions; the example above shows that a lower bound is not possible, in contrast to what happens with a single eigenfunction. \\

We deduce Theorem \ref{main} from a growth estimate for positive solutions of some higher order elliptic PDE. This estimate might be of independent interest and we state it as a theorem below.

\begin{theorem}
\label{doubling}
Let $m\geq 1$ be an integer. Consider the operator 
\begin{align}
\label{op}
    L=\sum_{i,j=1}^d a_{ij}(x)\partial_{x_i}\partial_{x_j}+\sum_{i=1}^d b_i(x)\partial_{x_i}+c(x)
\end{align}
on some open domain $\Omega\subset\R^d$, where we assume that all the coefficients are real-valued, $\mathcal{C}^{2m}$-smooth and bounded and the matrix $A=(a_{ij})$ is symmetric and uniformly positive definite. Let $\mathscr{L}$ be the operator of order $2m$ 
\begin{align}
\label{eq:1}
    \mathscr{L}=\prod_{k=1}^{m}(L+\gamma_k),
\end{align}
where $\gamma_k\in \R$ and $|\gamma_k|\leq \gamma$ for each $k$. Then there exists some $r_0$ which depends on $m$ and $L$ and $c>2$ depending on $L$, such that for every $r\leq r_0/c$ and every $u$ such that $u\in L^1_{loc}(B(x_0,cr))$, $u\geq 0$ and $\mathscr{L}u\leq 0$ in $B(x_0,cr)$, the following estimate holds:
\begin{align}
\label{doub}
    \int_{B(x_0,2r)}u \leq C \int_{B(x_0,r)}u,
\end{align}
where the constant $C=C(m,\gamma,L)$ is independent of $u$ and $r$. The explicit dependence of $C$ on $m$ and $\gamma$ is as follows:
\begin{enumerate}[(a)]
    \item if $\gamma_k\geq 0$ for any $k$, $C(m,\gamma,L)$ is independent of $\gamma$ and 
    \begin{align*}
    C(m,L)= \begin{cases}
    C_0\exp(C_1m^{\frac{d+1}{2}})\quad \text{if} \quad d\geq 4.\\
    C_0\exp(C_1m\log^2m)\quad \text{if}\quad d=1,2,3.      
    \end{cases}
\end{align*}
    \item if no sign condition is given on $\gamma_k$,
    \begin{align*}
    C(m,\gamma,L)= \begin{cases}
    C_0\exp(C_1(m^{\frac{d+1}{2}}+\sqrt{m\gamma}r_0))\quad \text{if} \quad d\geq 4.\\
    C_0\exp(C_1(m\log^2m+\sqrt{m\gamma}r_0))\quad \text{if}\quad d=1,2,3.      
    \end{cases}
\end{align*}
\end{enumerate}
where $C_0$ and $C_1$ vary from line to line but depend only on $L$.
\end{theorem}
A couple of remarks are in order. 
\begin{remark}
\textit{(i).} The product form \eqref{eq:1} of the higher order operator is not really necessary for an estimate like \eqref{doub} to hold; we state it here in this form because it allows us to have a better control over the constants involved, which is useful in the application to eigenfunctions. In section 5 we will state a more general, albeit less precise, version of Theorem 2. \\
\textit{(ii).} We think of Theorem \ref{doubling} as a substitute for the Harnack inequality for higher order operators. The pointwise Harnack inequality does not hold even for the euclidean bilaplacian, as the example $u(x)=x_1^2$ shows; Theorem \ref{doubling} shows that in average it is still true that positive solutions of higher order elliptic PDEs do not grow much. 
\end{remark}

In the next section we will show how to prove Theorem \ref{main} using Theorem \ref{doubling}. Section 3 will then be devoted to the proof of Theorem \ref{doubling}. In section 4 we indicate an alternative proof of Theorem \ref{main} which does not rely on Theorem \ref{doubling} and gives a slightly more precise result if the number of eigenfunctions in the sum is low. In the last section we discuss a form of Theorem \ref{doubling} for more general higher order elliptic operators and prove similar estimates for positive solutions of other equations using the same method.

\section{Extension tricks}
We start by describing a proof of the aforementioned density of the zero set of a single eigenfunction $\phi_\ld$. Consider the function $$h(x,t)=\phi_{\ld}(x) e^{\sqrt{\ld}t}$$
on the manifold $M\times \R$; note that $h$ is a harmonic function on $M\times \R$ and its zero set is a cylinder over the zero set of $\phi_{\ld}$. Assuming that $\phi_{\ld}$ is, say, positive in some geodesic ball of radius $2r$, the Harnack inequality applies to $h$ to give 
$$\sup_{B_r\times[-r,r]}h\leq C\inf_{B_r\times[-r,r]}h.$$
But $\inf_{B_r\times[-r,r]}h\leq\inf_{B_r}\phi_{\ld}$ and $\sup_{B_r\times[-r,r]}h\geq e^{\sqrt{\ld}r}\inf_{B_r}\phi_{\ld}$,
so that one obtains the restriction $r\leq\frac{C}{\sqrt{\ld}}$, where $C$ depends on $(M,g)$.\\

Let now $f=\sum_{j=1}^{m}a_j\phi_{{\ld}_j}$, where we stress again that $\ld_1$ is not the first eigenvalue of the Laplace-Beltrami operator. Note that $f$ satisfies the equation 
\begin{align}
    \left[\prod_{j=1}^m(\Delta_g+\ld_j)\right]f=0.
\end{align}
This observation was used in \cite{D}; we now combine it with an extension trick akin to the one mentioned above. We define the function 
\begin{align}
\label{extension}
h(x,t)=f(x)e^{\sqrt{\ld_1}t}
\end{align}
on $M\times \R$, so that the zero set of $h$ is a cylinder over the zero set of $f$. An elementary manipulation shows that $h$ satisfies:
\begin{align}
\label{eq main}
\Delta_{M\times \R}\left[\prod_{k=2}^{m}(\Delta_{M\times \R}+(\ld_k-\ld_1))\right]h=0.
\end{align}

Note that the Harnack inequality is not available as in the second order case, and that is why we need something like Theorem \ref{doubling} to handle equation \eqref{eq main}.\\

Pick a point $p\in M$ where $f(p)>0$. Consider a geodesic ball centered at $p$ of radius $R$, where we take $R$ smaller than the injectivity radius of $M$; take normal coordinates on this ball: if $R$ is small enough (depending on $g$) the geodesic and euclidean distances are comparable, so we will consider equivalent geodesic and euclidean balls. Recall that in local coordinates 
$$\Delta_g=\frac{1}{\sqrt{\det g}}\diver(\sqrt{\det g}g^{-1}\nabla(\cdot))$$ where divergence and gradient are in the euclidean sense, so that $\Delta_g=\sum_{i,j=1}^d g^{ij}(x)\partial_{x_i}\partial_{x_j}+(\det g)^{-\frac{1}{2}}\sum_{j=1}^d\sum_{i=1}^d\partial_{x_i}(\sqrt{\det g}g^{ij}(x))\partial_{x_j}$ is of the form \eqref{op} considered in Theorem \ref{doubling}.\\

Suppose now $f$ is positive on a ball $B(p,cr)$, $cr\leq R$; then $h$ is positive in $B(p,cr)\times [-cr,cr]$. We can assume that $r$ is small enough to apply the results of Theorem \ref{doubling}. Note that by hypothesis $\ld_k-\ld_1\geq 0$, so that we may use the better estimate for the constant in case \textit{(a)} of Theorem \ref{doubling}. We apply Theorem \ref{doubling} on $h$ (inscribing the cylinder in a ball of twice the radius and applying the doubling estimate twice) and obtain the inequality:
\begin{align*}
\int_{B(p,2r)\times [-2r,2r]}h \leq C_0\exp(C_1\sigma(m))\int_{B(p,r)\times [-r,r]}h,
\end{align*}
where the constants $C_0$ and $C_1$ depend on $(M,g)$, $\sigma(m)=m^{\frac{d+2}{2}}$ if $d\geq 3$ and $\sigma(m)=m\log^2m$ if $d=1,2$ as per Theorem 2. Recalling the definition of $h$, one obtains:
\begin{align*}
\int_{B(p,2r)}f\int_{-2r}^{2r}e^{\sqrt{\ld_1}t} \leq C_0\exp(C_1\sigma(m))\int_{B(p,r)}f\int_{-r}^{r}e^{\sqrt{\ld_1}t} \\
\leq C_0\exp(C_1\sigma(m))\int_{B(p,2r)}f\int_{-r}^{r}e^{\sqrt{\ld_1}t},
\end{align*}
and therefore
\begin{align*}
\exp(\sqrt{\ld_1}r)\leq C_0\exp(C_1\sigma(m)).
\end{align*}
This last inequality can only hold if $\sqrt{\ld_1}r\leq C\sigma(m)$. Since for a general smooth metric Theorem \ref{doubling} requires $r_0$ to depend on $m$, the above is only meaningful when 
$$\frac{C\sigma(m)}{\sqrt{\ld_1}}<r_0(m),$$
whence the condition $\ld_1>C_1(M,g,m)$. The proof of Theorem \ref{main} is thus complete.

\section{Proof of Theorem 2}
In order to simplify notation, we set $x_0=0$. We set $$\mathscr{L}=\prod_{k=1}^{m}(L+\gamma_k),$$ with $L=\sum_{i,j=1}^d a_{ij}(x)\partial_{x_i}\partial_{x_j}+\sum_{i=1}^d b_i(x)\partial_{x_i}+c(x)$ a uniformly elliptic second order operator with $C^{2m}$-smooth coefficients as in \eqref{op}. We denote by $A$ the matrix $(a_{ij})$. Recall the notation $\gamma=\max_k|\gamma_k|$. The strategy of the proof is to integrate by parts on $B_{3r}$ with a weight that has good ``convexity" properties under the formal adjoint $\mathscr{L}^*$. This is made precise in the following lemma, whose proof provides an explicit construction of the weight.

\begin{lemma}
\label{weight}
Suppose that $L$ is of the form \eqref{op}, with $A(0)=I$. Given $r\leq r_0$, where $r_0$ depends on $m$ and $L$, there exists $v_r\geq 0\in \mathcal{C}^{2m}(\overline{B_{3r}})$ such that:
\begin{enumerate}[(i)]
\item $\partial^{\mu}v_r=0$ on $\partial B_{3r}$ for $|\mu|\leq 2m-1$;
\item $\mathscr{L}v_r\geq 1$ on $B_{2r}\setminus B_r$;
\item $\mathscr{L}v_r\geq 0$ on $B_{3r}\setminus B_r$;
\item $\sup_{B_r}|\mathscr{L}v_r|\leq C(m,\gamma,L)$,
\end{enumerate}
where $C(m,\gamma,L)$ has the same dependency on $m$ and $\gamma$ as in the statement of Theorem \ref{doubling}, conditional on the sign of the $\gamma_k$'s.
\end{lemma}

\begin{remark}
We point out that the weight $v_r$ we construct depends on $r$, but the sup estimate in \textit{(iv)} does not. This is important, as the constant in Theorem \ref{doubling} does not depend on $r$.
\end{remark}

We postpone for a moment the proof of Lemma \ref{weight} and show how it provides the main ingredient for the proof of Theorem \ref{doubling}. By assumption $u$ is positive in a ball $B_{cr}$, where we choose $c$ so that $B_{cr}$ contains the interior of the ellipsoid defined by $\{A(0)^{-1}x\cdot x=3r\}$. We first perform a linear change of coordinates to get $A(0)=I$. Note that the top order part of $L^*$ is still $\sum_{i,j=1}^d a_{ij}(x)\partial_{x_i}\partial_{x_j}$, so that $L^*$ is of the form \eqref{op} and Lemma \ref{weight} applies to $\mathscr{L}^*=\Pi_{k=1}^{m}(L^*+\gamma_k)$. Using $\mathscr{L}u\leq 0$, $v_r\geq 0$ and property \textit{(i)} above, which eliminates the boundary terms, we integrate by part and have in the new coordinates:
\begin{align*}
0\geq\int_{B_{3r}}(\mathscr{L} u)v_r=\int_{B_{3r}}u(\mathscr{L}^* v_r).
\end{align*}
Combining this with \textit{(ii)}, \textit{(iii)} and \textit{(iv)}, and using the positivity of $u$, we estimate:
\begin{align*}
\int_{B_{2r}\setminus B_r}u\leq \int_{B_{2r}\setminus B_r}u(\mathscr{L}^*v_r)\leq \int_{B_{3r}\setminus B_r}u(\mathscr{L}^*v_r)\\
\leq-\int_{B_r}u(\mathscr{L}^*v_r)\leq C(m,\gamma,L) \int_{B_r}u.
\end{align*}
We thus obtained that 
\begin{align*}
\int_{B_{2r}}u=\int_{B_r}u+\int_{B_{2r}\setminus B_r}u\leq (1+C(m,\gamma,L))\int_{B_r}u
\end{align*}
in the new coordinates. Trasforming back to the original coordinates we get the estimate above for concentric ellipsoids of eccentricity governed by the ellipticity of $A$. Applying the inequality over ellipsoids several times (depending on the ellipticity) in order to be able to inscribe a ball in the largest ellipsoid and circumscribe a ball of half the radius to the smallest ellipsoid, we get the inequality over balls in the statement of Theorem \ref{doubling}. 

\subsection*{Construction of the weight}
The remainder of the present section is dedicated to the proof of Lemma \ref{weight}. We consider as a starting point the function
\begin{align}
    v(x)=\frac{1}{|x|^{\alpha}}.
\end{align}
In \cite{L} Landis shows that $Lv\geq 0$ in $B_{r_0}\setminus \{0\}$, where $r_0$ depends on $L$, if $\alpha$ is sufficiently large depending on the ellipticity constants of $L$ (as an example, for the euclidean laplacian $\alpha\geq d-2$ suffices). We now show that this positivity property is still true when we act upon $v$ by higher order operators such as $\mathscr{L}$. Since by assumption $A(0)=I$, we can write $L=\Delta_{eucl}+R$ in $B_{r_0}$, where $R$ is a second order operator whose top order coefficients are bounded by $C|x|$. \\

We define $\nu=x/|x|$ (it is the normal vector at the boundary of the ball), and note that $\nabla v=-\alpha |x|^{-(\alpha +1)}\nu$. Then 
\begin{align}
    Lv=\diver(\nabla v)+Rv=\frac{\alpha (\alpha +2-d)}{|x|^{\alpha +2}}+O\left(r_0 \frac{\alpha^2}{{|x|^{\alpha +2}}}\right).
\end{align}
If $r_0$ is small enough and $\alpha$ is large enough (depending on $L$), we get that $$Lv(x)\sim \frac{{\alpha}^2}{|x|^{\alpha+2}}$$ in $B_{r_0}\setminus \{0\}$ where here and in the sequel the symbol $\sim$ (as well as $\lesssim$ and $\gtrsim$) hides  constants that depend on $L$ only. Going further, a simple calculation shows that 
\begin{align*}
L^m v(x)=\frac{\alpha (\alpha+2-d)\dots(\alpha+2(m-1))(\alpha+2m-d)}{|x|^{\alpha+2m}}+R_mv,
\end{align*}
where $R_mv$ is a sum of terms all of which have at least a factor of $|x|$ more than the main term. By taking $r_0$ small enough (depending on $m$ as well as on $L$) and assuming $\alpha \gtrsim m$, we get that 
$$L^m v(x)\sim \frac{\prod_{j=0}^{2m-1}(\alpha+j)}{|x|^{\alpha+2m}}$$
in $B_{r_0}\setminus \{0\}$. We note that in general $r_0$ depends on $m$ in an unspecified way, because terms containing up to $2m-1$ derivatives of $A$ appear in $R_mv$.\\

The function $v$ is then a good candidate for our weight, but we need to modify it in order to deal with the derivatives at $\partial B_{3r}$ (as well as around $0$ where $v$ blows up). To this end, consider the Taylor polynomial centered at $s$ of order $l$ of the function $t^{-\alpha}$ for $t>0$ , call it $P_{l}(t;s)=\sum_{k=0}^l v_k(t;s)$. Note that the term $v_k$ of order $k$ is 
$$v_k(t;s)=\frac{\prod_{j=0}^{k-1}(\alpha+j)}{s^{\alpha+k}}\frac{(s-t)^k}{k!}$$
($v_0(t;s)=s^{-\alpha}$) which is positive for every $k\geq 0$ when $t\leq s$, so that $0\leq t^{-\alpha}-P_{l}(t;s)\leq t^{-\alpha}$ for any $l\geq 0$. We now set
\begin{align}
\widetilde{v_r}(x)=v(x)-P_{2m-1}(|x|;3r).
\end{align}
The following lemma is an easy consequence of the construction.
\begin{lemma}
$\widetilde{v_r}$ satisfies condition \textit{(i)} of Lemma 1, that is $\partial^{\mu}\widetilde{v_r}=0$ on $\partial B_{3r}$ for $|\mu| \leq 2m-1$.
\begin{proof}
Denoting by $\partial_{\nu}$ the normal derivative, we have that by construction $\partial_{\nu}^k \widetilde{v_r}=0$ on $\partial B_{3r}$ for $0\leq k \leq 2m-1.$ Now note that since $\widetilde{v_r}{\big|}_{\partial B_{3r}}=0$, the gradient at the boundary is parallel to the normal, i.e. $\nabla \widetilde{v_r}=(\partial_{\nu} \widetilde{v_r})\nu$ on $\partial B_{3r}$ and therefore all of the first derivative of $\widetilde{v_r}$ are zero at the boundary. The same argument holds for the higher derivatives as well: the vanishing of $\partial_{\nu}^{k-1} \widetilde{v_r}$ and $\partial_{\nu}^k \widetilde{v_r}$ at the boundary implies the vanishing of all the derivatives of order $k$ at the boundary, and the lemma is proved.
\end{proof}
\end{lemma}

We now need to check that the positivity property of $v$ under iterations of $L$ is preserved by the perturbation we introduced; this is clear if $\alpha$ is allowed to be very large with respect to $m$. In the next lemma we try to get some control over the minimal growth of $\alpha$ in $m$; if a reader is not interested in such things, it can be safely skipped.\\

Let us first note an elementary estimate: if $0\leq l\leq 2m-1$ is an integer, 
\begin{align}
    |\Delta_{eucl}^m|x|^l|\lesssim l!!(l+d-2)!!(2m-1-l)!|x|^{l-2m},
\end{align}
where we indicate with $k!!$ the semifactorial, that is the product from $1$ to $k$ of all the integers with the same parity as $k$.
\begin{lemma}
If $\alpha\gtrsim m$, the inequality 
\begin{align}
    |L^m P_{2m-1}(|x|;3r)|\lesssim \frac{m^{\frac{d+1}{2}}\prod_{j=0}^{2m-2}(\alpha+j)}{(3r)^{\alpha}|x|^{2m}}
\end{align}
holds for $x\in B_{3r}\setminus \{0\}$. In addition, if $d\in \{1,2,3\}$,
\begin{align}
   L^m P_{2m-1}(|x|;3r)\lesssim \frac{C(\alpha,m)}{(3r)^{\alpha}|x|^{2m-1}}. 
\end{align}
\begin{proof}
    Recall that $P_{2m-1}(|x|;3r)=\sum_{k=0}^{2m-1}\frac{\prod_{j=0}^{k-1}(\alpha+j)}{k!(3r)^{\alpha+k}}(3r-|x|)^k$. We write as before $L^m=\Delta_{eucl}^m+R_m$. Using (11) and expanding $(3r-|x|)^k=\sum_{l=0}^k(-1)^l\binom{k}{l}|x|^l(3r)^{k-l}$, we estimate:
    \begin{align*}
        |L^m P_{2m-1}(|x|;3r)|\lesssim\sum_{k=1}^{2m-1}\frac{\prod_{j=0}^{k-1}(\alpha+j)}{k!(3r)^{\alpha+k}}\times\\\times\sum_{l=1}^k\binom{k}{l}(3r)^{k-l}l!!(l+d-2)!!(2m-1-l)!|x|^{l-2m}\\
        \lesssim \frac{1}{(3r)^{\alpha}|x|^{2m}}\sum_{k=1}^{2m-1}c_k\prod_{j=0}^{k-1}(\alpha+j),
    \end{align*}
    where $c_k=\sum_{l=1}^k\frac{(l+d-2)!!}{(l-1)!!}\frac{(2m-1-l)!}{(k-l)!}$. Observe now that $$\frac{(l+d-2)!!}{(l-1)!!}\lesssim (l+d-2)^{\frac{d-1}{2}}.$$ We write $$\frac{(2m-1-l)!}{(k-l)!}=(2m-1-k)!\binom{2m-1-l}{2m-1-k}$$ and note that we have the identity  $$\sum_{l=1}^k\binom{2m-1-l}{2m-1-k}=\binom{2m-1}{2m-k},$$
    so that $c_k\lesssim m^{\frac{d-1}{2}}\frac{(2m-1)!}{(2m-k)(k-1)!}$. We can then estimate:
    \begin{align*}
        \sum_{k=1}^{2m-1}c_k\prod_{j=0}^{k-1}(\alpha+j)\lesssim m^{\frac{d-1}{2}}\sum_{k=1}^{2m-1}\frac{(2m-1)!}{(k-1)!}\prod_{j=0}^{k-1}(\alpha+j)\\ \lesssim
        m^{\frac{d-1}{2}}(2m-1)!\frac{\prod_{j=0}^{2m-2}(\alpha+j)}{(2m-2)!}
        \lesssim m^{\frac{d+1}{2}}{\prod_{j=0}^{2m-2}(\alpha+j)},
    \end{align*}
    where the second inequality holds because $\sum_{k=1}^{2m-1}\frac{\prod_{j=0}^{k-1}(\alpha+j)}{(k-1)!}$ can be estimated by, say, twice its last term if $\alpha\gtrsim m$, and (12) is proved.
    
    Consider now $d\in \{1,2,3\}$: for $0\leq l\leq 2m-1$, $(-1)^l\Delta_{eucl}^m|x|^l\leq 0$ in $B_{3r}\setminus B_r$ (if $l$ is even, equality holds trivially in any dimension; when $d=1$ equality holds; when $d=3$ the inequality holds because $|x|^{-1}$ is the fundamental solution for the laplacian; when $d=2$, the inequality holds because $\Delta_{eucl}|x|^{-1}=|x|^{-3}\geq 0$). Again by expanding $(3r-|x|)^k=\sum_{l=0}^k(-1)^l\binom{k}{l}|x|^l(3r)^{k-l}$, this implies that $\Delta_{eucl}^m P_{2m-1}(|x|;3r)\leq 0$, so that $L^m P_{2m-1}(|x|;3r)\leq R_m P_{2m-1}(|x|;3r)$, from which (13) follows.
\end{proof}
\end{lemma}

Recall now that $$L^m v(x)\gtrsim \frac{\prod_{j=0}^{2m-1}(\alpha+j)}{|x|^{\alpha+2m}}$$ for $x\in B_{3r}\setminus B_r$. Since $\widetilde{v_r}=v-P_{2m-1}(|x|;3r)$, by Lemma 3
\begin{align*}
    L^m\widetilde{v_r}(x)\gtrsim \frac{\prod_{j=0}^{2m-1}(\alpha+j)}{|x|^{\alpha+2m}}-\frac{m^{\frac{d+1}{2}}\prod_{j=0}^{2m-2}(\alpha+j)}{(3r)^{\alpha}|x|^{2m}},
\end{align*}
so that if $\alpha\gtrsim m^{(d+1)/2}$
\begin{align}
\label{main term}
    L^m\widetilde{v_r}(x)\gtrsim \frac{\prod_{j=0}^{2m-1}(\alpha+j)}{|x|^{\alpha+2m}}
\end{align} 
for $x\in B_{3r}\setminus B_r$ (with a smaller implied constant). When $d\in \{1,2,3\}$, we can use the better estimate (13) so that (14) holds with $\alpha\gtrsim m$.\\

It is now time to consider the full operator $\mathscr{L}$ and see what changes from only its top order part $L^m$. We start by noting the following:

\begin{lemma}
Under the above constraints on $\alpha$ and $r_0$ and in the above notations, $L^q \widetilde{v_r}\geq 0$ on $B_{3r}\setminus \{0\}$ for any $0\leq q\leq m-1$. Moreover, $L^q \widetilde{v_r}\lesssim L^q v\lesssim (C\alpha)^{2q} |x|^{-(\alpha+2q)}$ on $B_{3r}\setminus \{0\}$.
\begin{proof}
Call for ease of notation $w_q=L^q \widetilde{v_r}$. Note that $w_q$ and all of its derivatives up to order $2(m-q)-1$ are zero on $\partial B_{3r}$, since they are combinations of derivatives up to order $2m-1$ of $\widetilde{v_r}$. Let $\partial_{\rho}$ indicate derivative in the radial direction; then the estimate (11) with the same right hand side holds for $\partial_{\rho}^{2(m-q)}\Delta_{eucl}^q |x|^l$ and then with the same proof of Lemma 3 one gets $$\partial_{\rho}^{2(m-q)}w_q>0$$ on $B_{3r}\setminus \{0\}$. Fix $\theta\in\partial B_{3r}$; the above implies that $\partial_{\rho}^{2(m-q)-1}w_q(\rho,\theta)$ is strictly increasing for $0<\rho\leq 3r$ and, being zero for $\rho=3r$, is thus negative for $0<\rho\leq 3r$. This in turn implies that $\partial_{\rho}^{2(m-q)-2}w_q(\rho,\theta)>0$ for $0<\rho\leq 3r$, and so on. After $2(m-q)$ steps, we get $w_q(\rho,\theta)\geq 0$ for every $\theta\in\partial B_{3r}$. 

For the second part, call $f_{\alpha,l}(t)=t^{-\alpha}-P_l(t;s)$ in the same notations as in the definition of $\widetilde{v_r}$; as noted above, $0\leq f_{\alpha,l}(t)\leq t^{-\alpha}$ for any $l\geq 0$, and an elementary calculation shows $$\frac{d}{dt} f_{\alpha,l}=-\alpha f_{\alpha+1,l-1}.$$
From this it follows that $$\left|\frac{d^j}{dt^j} f_{\alpha,l}\right|\leq \left|\frac{d^j}{dt^j}t^{-\alpha}\right|$$ for any $j$,
which implies the statement.
\end{proof}
\end{lemma}

\noindent\textit{Case a).} As a consequence of Lemma 4, when all the $\gamma_k$'s in the definition of $\mathscr{L}$ are all positive, we obtain that $$\mathscr{L}\widetilde{v_r}\geq 0$$
on $B_{3r}\setminus B_r$ with no additional constraints on $\alpha$. On the smaller annulus $B_{2r}\setminus B_r$ we have that $$L^q \widetilde{v_r}\gtrsim \frac{\alpha^{2q}}{|x|^{\alpha+2q}}$$ for $0\leq q\leq m-1$. This implies that 
\begin{align*}
    \mathscr{L}\widetilde{v_r}(x)\gtrsim \frac{\prod_{k=1}^{m}(\alpha^2+\gamma_k|x|^2)}{|x|^{\alpha+2m}}
\end{align*}
for $x\in B_{2r}\setminus B_r$. In particular,
\begin{align}
    \inf_{B_{2r}\setminus B_r}\mathscr{L}\widetilde{v_r}\gtrsim (2r)^{-(\alpha+2m)}\prod_{k=1}^{m}(\alpha^2+\gamma_k r^2).
    \end{align}

\noindent\textit{Case b).} When the $\gamma_k$'s are allowed to have any sign, \eqref{main term} and the second part of Lemma 4 give (recall that $\gamma=\max|\gamma_k|$):
\begin{align*}
    \mathscr{L}\widetilde{v_r}(x)\gtrsim \frac{2\alpha^{2m}-(\alpha^2+\gamma|x|^2)^m}{|x|^{\alpha+2m}}
\end{align*}
for $x\in B_{3r}\setminus B_r$. We want $2\alpha^{2m}-(\alpha^2+\gamma|x|^2)^m\gtrsim (\alpha^2+\gamma r_0^2)^m$, which requires the additional condition $\alpha\gtrsim \sqrt{m\gamma}r_0$. In this case one obtains
\begin{align*}
    \mathscr{L}\widetilde{v_r}(x)\gtrsim \frac{(\alpha^2+\gamma r_0^2)^{m}}{|x|^{\alpha+2m}}
\end{align*}
for $x\in B_{3r}\setminus B_r$. In particular we have that
\begin{align}
   \inf_{B_{2r}\setminus B_r}\mathscr{L}\widetilde{v_r}\gtrsim (2r)^{-(\alpha+2m)}(\alpha^2+\gamma r_0^2)^{m}.  
\end{align}\\

The only task remaining is to define the weight inside $B_r$. To this end, consider a $\mathcal{C}^{\infty}$ radial cut-off function $\psi_r$, such that $0\leq \psi_r\leq 1$, $\psi_r\equiv 0$ in $B_{r/2}$ and $\psi_r\equiv 1$ in $(B_r)^c$. In H{\"o}rmander's treatise \cite{H}, Chapter 1, it is shown how to construct such a function in a way that 
$$|\partial^{\beta}\psi_r|\leq (C|\beta|\log^2|\beta|)^{|\beta|}r^{-|\beta|},$$ where $C$ is a dimensional constant; briefly, this is achieved through repeated convolutions of characteristic functions. Recall that for $\widetilde{v_r}$ we have $|\partial_{\rho}^k\widetilde{v_r}|\leq (C\alpha)^k|x|^{-(\alpha+k)}$; we then obtain for $q\leq m$
\begin{align*}
    |L^q(\psi_r\widetilde{v_r})|\lesssim \frac{C^{2q}(q\log^2 q+\alpha)^{2q}}{r^{\alpha+2q}}
\end{align*}
in $B_r$. We want the numerator to be controlled by $(C\alpha)^{2q}$ for every $q\leq m$, which requires the condition $\alpha\gtrsim m\log^2m$. This is already ensured when $d\geq 4$ by the constraint $\alpha\gtrsim m^{(d+1)/2}$ we set before; when $d\in \{1,2,3\}$, we previously required only $\alpha\gtrsim m$, so we need to assume now $\alpha\gtrsim m\log^2m$. With these requirements, we obtain 
\begin{align}
    \sup_{B_r}|\mathcal{L}(\psi_r\widetilde{v_r})|\lesssim \frac{C^m\Pi_{k=1}^m(\alpha^2+|\gamma_k| r^2)}{r^{\alpha+2m}}.
\end{align}
\\
We are finally ready to define 
\begin{align}
    v_r=\left(\inf_{B_{2r}\setminus B_r}\mathscr{L}\widetilde{v_r}\right)^{-1}\psi_r\widetilde{v_r}.
\end{align}
With $v_r$ defined as such, points \textit{(i)}, \textit{(ii)} and \textit{(iii)} in Lemma \ref{weight} are clear, and for point \textit{(iv)} it is enough to note that by (17) and (15) or (16) respectively
\begin{align*}
\sup_{B_r}|\mathscr{L}v_r|\lesssim\frac{C^m(2r)^{\alpha+2m}}{r^{\alpha+2m}}\lesssim\exp{C_1\alpha}.
\end{align*}
By choosing 
\begin{align*}
    \alpha\sim\begin{cases}
    m^{\frac{d+1}{2}}\quad \text{if} \quad d\geq 4.\\
    m\log^2m\quad \text{if}\quad d=1,2,3.
    \end{cases}
\end{align*}
when $\gamma_k\geq 0$ for every $k$ and 
\begin{align*}
    \alpha\sim\begin{cases}
    m^{\frac{d+1}{2}}+\sqrt{m\gamma}r_0\quad \text{if} \quad d\geq 4.\\
    m\log^2m+\sqrt{m\gamma}r_0\quad \text{if}\quad d=1,2,3.
    \end{cases}
\end{align*}
otherwise, we obtain a constant $C(m,\gamma,L)$ which has the dependencies stated in Theorem \ref{doubling}, and the proof of Lemma \ref{weight} is complete.

\section{A second proof of the density of zero sets}
We now present an alternative proof of a variant of Theorem 1 that does not rely on the extension trick or on the growth estimate contained in Theorem 2. For simplicity we state and prove the theorem for two eigenfunctions only. 
\begin{theorem}
    Let $\Omega\in \R^d$ be an open set and let $L$ be a uniformly elliptic second-order operator with smooth coefficients as in \eqref{op}. Consider any solution of the inequality $(L+\ld_1)(L+\ld_2)f\leq 0$ in $\Omega$. Then there exists a constant $C(L)$ such that, if $\ld_1$ and $\ld_2$ are large enough, $f$ has a zero in every ball of radius $$r^*=C(L)\sqrt{\frac{1}{\ld_1}+\frac{1}{\ld_2}}.$$
    \begin{proof}
    We argue by contradiction: suppose that $f$ is positive in a ball $B_r$ of radius $r>r^*$. We can assume $0\in \Omega$ and  that $B_r$ is centered at $0$. Consider the function 
    \begin{align}
        w(x)=(|x|^2-r^2)^k
    \end{align}
    for $x\in B_r$ and $k>4$ an even integer. Note that $\partial_{\nu}^j w{\big|}_{\partial B_r}=0$ for $0\leq j\leq k-1$; by the same argument as in Lemma 2 we then have 
    \begin{align}
        \partial^{\mu}w{\big|}_{\partial B_r}=0 
    \end{align}
    for $0\leq |\mu|\leq k-1$. We assume $r$ small enough and by a linear transformation of coordinates $L^*=\Delta_{eucl}+R$, where $R$ is a second-order operator whose top order coefficients are bounded by $C|x|$. We can then estimate:
    \begin{align*}
        L^*w(x)=2k(|x|^2-r^2)^{k-2}\{(k-1)|x|^2+d(|x|^2-r^2)\}+\\O(r^3(|x|^2-r^2)^{k-2}),
    \end{align*}
    and 
    \begin{align*}
     (L^*)^2w(x)=4k(k-1)(|x|^2-r^2)^{k-4}\{4(k-2)(k-3)|x|^4+\\4(k-2)(d+2)(|x|^2-r^2)|x|^2+d(d+2)(|x|^2-r^2)^2\}+\\O(r^5(|x|^2-r^2)^{k-4}).  
    \end{align*}
    Note that when, say, $|x|>r/2$, $L^*w>0$ and $(L^*)^2w>0$ if $k$ is large enough depending on $L$. This implies that $$(L^*+\ld_1)(L^*+\ld_2)w>0$$ in $B_r\setminus B_{r/2}$. We claim that if $r>r^*$ then $(L^*+\ld_1)(L^*+\ld_2)w>0$ also in $B_{r/2}$. In fact, using the formulas above, when $x\in B_{r/2}$ we have 
    \begin{align}
        (L^*+\ld_1)(L^*+\ld_2)w\gtrsim r^{2k-4}\{-C_0 k^3-C_1(\ld_1+\ld_2)kr^2+C_2\ld_1 \ld_2 r^4\}+\\
        O((\ld_1+\ld_2)r^{2k-1})+O(r^{2k-3}). \nonumber
    \end{align}
    Now, our assumption $r>r^*$ implies that $\ld_1\ld_2 r^2>C(\ld_1+\ld_2)$ for a sufficiently large $C$ depending on $L$; since $k$ depends on $L$ only, we conclude that the right hand side of (21) is positive. We thus found $w$ such that $(L^*+\ld_1)(L^*+\ld_2)w>0$ on $B_r$.
    Using (20), we obtain:
    \begin{align*}
        0\geq\int_{B_r}[(L+\ld_1)(L+\ld_2)f] w=\int_{B_r}f[(L^*+\ld_1)(L^*+\ld_2)w],
    \end{align*}
    and we reach a contradiction if $f$ is positive in $B_r$.
    \end{proof}
\end{theorem}
\begin{remark}
Note that the technique above can also be used to give another proof of the density of zero sets for a single eigenfunction. Both the proof above and the one illustrated before can easily be modified to show density of zero sets for eigenfunctions (and linear combinations thereof) of operators of the form $\Delta^{2q+1}$, where $q$ is a positive integer. Neither of the proofs work for $\Delta^2$, or any other operator with positive eigenvalues. 
\end{remark}

\section{Other higher-order operators}
We stated Theorem \ref{doubling} for a higher-order operator which was the product of second-order operators, because in that case we have better control of the constants needed for the application to linear combinations of eigenfunctions. The product form is however not necessary for the proof to go through and a doubling estimate can be proved for positive (sub)-solutions of higher-order elliptic operators in other forms. We now state a theorem in this direction and give a sketch of the proof that follows the argument of Section 3, in some sense more easily since we do not keep track of the rate of growth of the constants. A somewhat similar result is contained in \cite{KKM}, Theorem 3 there, where non-negative solutions of higher-order elliptic PDEs in a cone are considered. \\

Let $\mathscr{L}$ be a linear differential operator of order $2m$ in a domain $\Omega\in \R^d$. We write in a convenient way $$\mathscr{L}=\sum_{|\mu|=2m}a_{\mu}(x)\partial^{\mu}+\sum_{|\mu|<2m}b_{\mu}(x)\partial^{\mu},$$ so that we can write the formal adjoint in the form $$\mathscr{L}^*=\sum_{|\mu|\leq 2m}a_{\mu}(x)\partial^{\mu}.$$ We assume the coefficients satisfy the uniform ellipticity condition 
$$\sum_{|\mu|=2m}a_{\mu}(x)\xi^{\mu}\geq C_1|\xi|^{2m}$$ for some $C_1>0$ and any $\xi\in \R^d$ and a bound $|a_{\mu}(x)|\leq C_2$. 
\begin{theorem}
    Let $\mathscr{L}$ be as above, and consider a ball $B(x_0,4R)\subset \Omega$, $4R\leq 1$. For any $u\in L^1_{loc}(B(x_0,4R))$ such that $u\geq 0$ and $\mathscr{L}u\leq 0$ in $B(x_0,4R)$, the following estimate holds:
    \begin{align}
        \int_{B(x_0,2R)}u\leq C \int_{B(x_0,R)}u,
    \end{align}
    where the constant $C$ depends on $m$, $C_1$ and $C_2$.
  \begin{proof}[Proof (Sketch)]
  We suppose $x_0=0$. The proof goes in the same way as that of Theorem \ref{doubling}. We set $v(x)=|x|^{-\alpha}$, where $\alpha$ is a large positive number that will be chosen later. Take a multi-index $\mu$ with $|\mu|=2m$; an elementary calculation shows that 
  $$\partial^{\mu}v(x)=\frac{\prod_{j=0}^{2m-1}(\alpha+2j)}{|x|^{\alpha+2m}}\left(\frac{x}{|x|}\right)^{\mu}+g(x),$$
  where $g(x)$ has a lower order in $\alpha$ than the first term (and the same in $|x|$). As in the proof of Theorem \ref{doubling} we want to integrate by parts with weight $v$ on the ball $B_{3R}$, so we modify it by subtracting its Taylor polynomial at $|x|=3R$; that is, we define
  $$\widetilde{v_R}(x)=v(x)-\sum_{k=0}^{2m-1}\frac{\prod_{j=0}^{k-1}(\alpha+j)}{(3R)^{\alpha+k}}\frac{(3R-|x|)^k}{k!}.$$
  Note that each polynomial term only has up to $2m-1$ factors containing $\alpha$; we then have 
  $$\mathscr{L}^*\widetilde{v_R}=\frac{\prod_{j=0}^{2m-1}(\alpha+2j)}{|x|^{\alpha+2m}}\sum_{|\mu|=2m}a_{\mu}(x)\left(\frac{x}{|x|}\right)^{\mu}+g_{\alpha}(x),$$
  where again $g_{\alpha}$ is lower order in $\alpha$. By ellipticity $$\sum_{|\mu|=2m}a_{\mu}(x)\left(\frac{x}{|x|}\right)^{\mu}\geq C_1,$$ so that if we choose $\alpha$ large enough (depending on $m$, $C_1$ and $C_2$) we obtain 
  $$\mathscr{L}^*\widetilde{v_R}\gtrsim \frac{\prod_{j=0}^{2m-1}(\alpha+2j)}{|x|^{\alpha+2m}}$$
  in $B_{3R}\setminus\{0\}$. We now have to modify $\widetilde{v_R}$ close to $0$ so that it is $\mathcal{C}^{2m}$ in $B_{3R}$; we can do that for instance by choosing a cut-off function $\psi$ which is equal to $0$ in $B_{R/2}$ and equal to $1$ in $B_{3R}\setminus B_R$, and defining $v_R=\psi \widetilde{v_R}$. Using the hypothesis of the theorem we estimate 
  \begin{align*}
      0\geq \int_{B_{3R}}(\mathscr{L}u) v_R=\int_{B_{3R}}u\mathscr{L}^*v_R\geq \int_{B_{R}}u\mathscr{L}^*v_R+\int_{B_{2R}\setminus B_R}u\mathscr{L}^*v_R\\
      \geq -C_{\alpha}R^{-(\alpha+2m)}\int_{B_{R}}u+C^{'}_{\alpha}(2R)^{-(\alpha+2m)}\int_{B_{2R}\setminus B_R}u,
  \end{align*}
  so that $\int_{B_{2R}\setminus B_R}u\leq C^{''}_{\alpha}\int_{B_R}u$, from which (22) follows.
  \end{proof}  
\end{theorem}

It is clear from the proof that it is not necessary for the lower order coefficients to be bounded, and it would suffice to have $|a_{\mu}(x)|\leq C_2|x|^{|\mu|-2m}$, a condition which appears in \cite{KKM} on cones. We illustrate a similar situation with an example where we have better control of the constants.

\subsection*{Positive solutions of a Schr{\"o}dinger equation}
Consider the equation 
\begin{align}
\label{sc}
    Lu+Vu=0
\end{align}
in $B_R\subset \R^d$, $d\geq 3$, where $L$ is a uniformly elliptic second order operator as in \eqref{op}, and suppose that there exists some $\eta>0$ such that the potential $V$ satisfies
\begin{align}
\label{pot}
    |V(x)|\leq \frac{\eta}{|x|^2}.
\end{align}
Equation (23) with assumptions on $V$ like (24) was studied in \cite{GL} within the context of unique continuation. We can recover a result of \cite{CFG} (Lemma 2.3 there), which is used as an intermediate step towards the full Harnack inequality for positive solutions of \eqref{sc}. In \cite{CFG}, $V$ is assumed to belong to the so-called Stummel class, which means that 
\begin{align}
    \lim_{r\to 0}\sup_{x\in\Omega}\int_{|x-y|<r}\frac{|V(y)|}{|x-y|^{d-2}}dy=0
\end{align}
for any bounded set $\Omega\in \R^d$.
The assumption \eqref{pot} allows a faster growth of $V$ near one singularity; in fact, $u(x)=|x|^2$ is a solution of (23) with potential satisfying (24) but not (25), and the full Harnack inequality does not hold for $u$. Our proof is quite different from the one of \cite{CFG} and follows the methods of the preceding sections.

\begin{proposition}
Let $u\in W^{2,\infty}_{loc}(B_R)$ be a strong solution of \eqref{sc}, with $V$ satisfying \eqref{pot}. Suppose that $u\geq0$ in a ball of radius $4r\leq r_0$, where $r_0$ depends on $L$. Then 
\begin{align}
    \int_{B_{2r}}u \leq \exp{(C_0+C_1\sqrt{\eta})}\int_{B_r}u
\end{align}
where $C_0$ and $C_1$ depend on $L$ only.
\begin{proof} The proof runs in the same way as that of Theorem \ref{doubling}, more easily since the operator is just of second order: integrate by parts with weight $v_r$ such that $v_r$ is a supersolution for $L^*+V$ on the annulus $B_{2r}\setminus B_r$. We can use $v_r$ constructed as in Lemma \ref{weight}; explicitly, let
$$\widetilde{v_r}(x)=\frac{1}{|x|^{\alpha}}-\frac{1}{(3r)^{\alpha}}+\frac{\alpha}{(3r)^{\alpha+1}}(|x|-3r).$$
We already estimated $L^*\widetilde{v_r}\gtrsim \alpha^2|x|^{-(\alpha+2)}$ if $\alpha$ is large enough depending on $L$, and noted $0\leq \widetilde{v_r}(x)\leq |x|^{-\alpha}$, so that by \eqref{pot}
\begin{align}
    (L^*+V)\widetilde{v_r}\gtrsim \frac{\alpha^2}{|x|^{\alpha+2}}-\frac{|V(x)|}{|x|^{\alpha}} \geq \frac{\alpha^2-\eta}{|x|^{\alpha+2}}
\end{align}
in $B_{2r}\setminus B_r$. Pick $\alpha\sim C_0+C_1\sqrt{\eta}$, with $C_0$ and $C_1$ large enough. As in the conclusion of the proof of Lemma 1 we can modify $\widetilde{v_r}$ inside $B_r$ by multiplying with a smooth cut-off function $\psi_r$ which is identically $0$ in $B_{r/2}$ and $1$ outside $B_r$. Let $v_r=\psi_r\widetilde{v_r}$; we have
\begin{align}
    \sup_{B_r}|(L^*+V)v_r|\lesssim \frac{\alpha^2+\eta}{|x|^{\alpha+2}}.
\end{align}
Recall that the constant in the integral estimate is given by the quotient of an upper bound on $B_r$ and a lower bound on $B_{2r}\setminus B_r$ for $(L^*+V)v_r$; in our case, using (27) and (28), this is $\sim 2^{\alpha+2}$. The statement then follows by the choice of $\alpha$.
\end{proof}
\end{proposition}
\begin{remark}
Note that the function $u(x)=|x|^{2k}$, where $k$ is a positive integer, shows that the order of growth of the constant in (26) with respect to $\eta$ is sharp.
\end{remark}

\section*{Acknowledgements}
The author warmly thanks Alexander Logunov for suggesting the problem and for helpful comments on previous drafts. He also gratefully acknowledges the continuous support and guidance of Eugenia Malinnikova, whose many comments and suggestions greatly helped improve the present work. 

This work was started while the author was visiting the Department of Mathematics at Princeton University, and parts of it were completed while he was a Visiting Student Researcher at the Department of Mathematics at Stanford University; it is a pleasure to thank both institutions for their hospitality.

Thanks are also given to the anonymous referee for helpful comments on the presentation and on the history and bibliography of some facts discussed in the introduction. 

The author is supported by Project 275113 of the Research Council of Norway.

\end{document}